\documentstyle[amssymb,12pt]{amsart}

\newcommand{\Chow}{\mbox{\it Chow}\,}
\newcommand{\dotc}{\setlength{\unitlength}{1pt}\begin{picture}(3,6)(1,1)
	\put(2,3){\circle*{2}}\end{picture}}
\newcommand{\DOT}{\setlength{\unitlength}{1pt}\begin{picture}(2.5,2)(1,1)
	\put(1,2){\circle*{1.7}}\end{picture}}
\newcommand{\Fdot}{{F\!_{\DOT}}}

\newcommand{\QED}{
\setlength{\unitlength}{1.0pt}%
\begin{picture}(7.5,10)
\put(0,-5){\rule{2.5pt}{5pt}}
\put(0,0){\rule{5pt}{2.5pt}}
\put(0,2.5){\rule{7.5pt}{2.5pt}}
\end{picture}\vspace{15pt}}

\newcommand{\DIAG}{\setlength{\unitlength}{3.5pt}\begin{picture}(6,6)(-6,-6)
\put( 0, 0){\line(0,-1){2}}
\put( 0,-4){\line(0,-1){2}}
\put(-1, 0){\line(0,-1){6}}
\put(-2, 0){\line(0,-1){6}}
\put(-3, 0){\line(0,-1){4}}
\put(-4,-2){\line(0,-1){4}}
\put(-5, 0){\line(0,-1){2}}
\put(-5,-4){\line(0,-1){2}}
\put(-6, 0){\line(0,-1){2}}

\put( 0, 0){\line(-1,0){1}}
\put(-2, 0){\line(-1,0){1}}
\put(-5, 0){\line(-1,0){1}}
\put( 0,-1){\line(-1,0){1}}
\put(-2,-1){\line(-1,0){1}}
\put(-5,-1){\line(-1,0){1}}
\put( 0,-2){\line(-1,0){4}}
\put(-5,-2){\line(-1,0){1}}
\put(-1,-3){\line(-1,0){3}}
\put( 0,-4){\line(-1,0){5}}
\put( 0,-5){\line(-1,0){2}}
\put(-4,-5){\line(-1,0){1}}
\put( 0,-6){\line(-1,0){2}}
\put(-4,-6){\line(-1,0){1}}
\end{picture}}

\newtheorem{lemma}{Lemma}
\newtheorem{thm}[lemma]{Theorem}

\newtheorem{cor}[lemma]{Corollary}

\begin{document}

\title[enumerative geometry for real varieties]{Enumerative geometry
for real varieties} 

\author{Frank Sottile}

\address{Mathematical Sciences Research Institute\\
	1000 Centennial Drive\\
	Berkeley, CA 94720\\
	USA}
\address{Department of Mathematics\\
        University of Toronto\\
        100 St. George St.\\
	Toronto, Ontario  M5S 3G3\\
	Canada (on leave)}
\email{sottile@@math.toronto.edu}
\date{16 July 1996}
\thanks{Research supported in part by NSERC grant \# OGP0170279}
 \thanks{This manuscript is based on author's talk at the 1995
AMS Summer Research Institute in Algebraic Geometry at Santa Cruz.
It will appear in the Proceedings and Symposia in Pure Mathematics
Series of the AMS.} 
\thanks{Extended version available by request from the author, or at
http://www.msri.org/members/bio/sottile.html}

\maketitle

\section{Introduction}
Of the geometric figures in a given family satisfying real
conditions, 
some figures are real while the rest occur in complex conjugate pairs,
and the  distribution of the two types depends subtly upon 
the configuration of the conditions.
Despite this difficulty, 
applications~(\cite{Byrnes},\cite{Mourrain_MEGA94},\cite{Ronga_Vust}) may
demand real solutions.
Fulton~\cite{Fulton_introduction_intersection} asked how many
solutions of an enumerative problem can be real, and 
we consider a special case of his question: 
Given a problem of enumerative geometry, are there real conditions 
such that every figure satisfying them is real?
Such an enumerative problem is {\em fully real}.

B\'{e}zout's Theorem, or rather the problem of intersecting
hypersurfaces in  ${\Bbb P}^n$,  is fully real.
This is readily seen for ${\Bbb P}^2$, and the argument generalizes 
to ${\Bbb P}^n$.
Suppose $X_0$ consists of $d$ real lines, $Y_0$ of $e$ real lines,
and $X_0$ meets $Y_0$ transversally in (necessarily) $d\cdot e$ 
real points.
Let $X$ and $Y$ be defined by suitably small
generic real deformations of the forms defining $X_0$ and $Y_0$.
Then  $X$ and $Y$   are smooth real plane
curves of degrees $d$ and $e$ meeting
transversally in $d\cdot e$ real points.

This argument used 
a degenerate case free of multiplicities;
$X_0$ and $Y_0$ are reduced and meet transversally.
While it is typical to introduce multiplicities 
(for example, in the proof of B{\'e}zout's Theorem in~\cite{Mumford_cpv})
to establish enumerative formulas,
multiplicities may lead to complex conjugate pairs
of solutions, complicating the search for real solutions.

All Schubert-type enumerative problems involving 
lines in ${\Bbb P}^n$ are fully real~\cite{sottile_real_lines}.
This follows from the existence of (multiplicity-free) deformations of 
generically transverse intersections of Schubert varieties into
sums of Schubert varieties.
Refining this method of multiplicity-free
deformations~\cite{sottile_mega96}  yields techniques for showing
other enumerative problems are fully real. 
Fulton, and more recently, Ronga, Tognoli, and
Vust~\cite{Ronga_Tognoli_Vust}, 
have shown the problem of 3264 conics tangent to five general 
plane conics is fully real.
Their analysis utilizes degenerate
conditions having multiplicities.

Enumerative problems that we know are not fully real share a
common flaw:  
they do not involve intersecting general subvarieties.
For example, Klein~\cite{Klein} showed that at most $n(n-2)$ of the
$3n(n-2)$ flexes on a real plane curve of degree $n$ can be real.
These flexes are the intersection of the curve with its Hessian
determinant,  {\em not} with a general curve of degree $3(n-2)$. 
Also, Khovanskii~\cite{Khovanskii_fewnomials} showed that if 
hypersurfaces in a complex torus are defined by polynomials 
with few monomials, then the real points of
intersection are at most a fraction of all points of
intersection. 
However, these are not generic hypersurfaces with given Newton polytope. 

Little is known about fully real enumerative problems.
For instance, we are unaware of a good theoretical framework for
studying fully real enumerative problems.
Also, it is not known how common it is for an enumerative problem to
be fully real.  
Here are some examples of enumerative problems worth considering:

\begin{enumerate}
\item Is the
Kouchnirenko-Bernstein Theorem~(\cite{Bernstein},~\cite{Kouchnirenko})
for hypersurfaces in a torus fully real?
That is, given lattice polytopes $\Delta_1,\ldots,\Delta_n$ in 
${\Bbb Z}^n$ are there real polynomials $f_1,\ldots,f_n$ where
$\Delta_i$ is the Newton polytope of $f_i$ and all solutions to the system
$f_1=\cdots=f_n=0$ in $({\Bbb C}^\times)^n$ are real?

\item Generalize the results
of~\cite{sottile_real_lines} and~\cite{sottile_mega96}: 
Are other (all?) Schubert-type enumerative problems on
flag varieties fully real?

\item All known examples involve spherical 
varieties~(\cite{Brion_spherical_introduction},%
~\cite{Knop_spherical_expository},~\cite{Luna_Vust_plongements}).
Which enumerative problems on other spherical varieties are fully real?

\item In~\cite{sottile_mega96} 
all problems of enumerating lines incident upon subvarieties of 
fixed dimension and degree in ${\Bbb P}^n$ are shown to be 
fully real. 
What is the situation for rational curves of higher degree?
(Degree 0 is B{\'e}zout's Theorem.)
For example, for which positive integers $d$ do there exist $3d-1$
real points in ${\Bbb P}^2$ such that the Kontsevich number $N_d$ of
degree $d$ rational curves passing through these
points~(\cite{Kontsevich_Manin}~\cite{Ruan_Tian}) are all real?
For an introduction to these questions of quantum cohomology, see the
paper by Fulton and Pandharipande~\cite{Fulton_Rahul} in this volume.

\end{enumerate}

This technique of multiplicity-free deformations may have applications
beyond showing the existence of real solutions.
When the deformations are explicitly described (which is the case in
most known results), it may be possible to obtain explicit solutions to
the enumerative problem using continuation methods of numerical
analysis~\cite{Allgower_Georg_1990} to follow real points in the
degenerate configuration backwards along the deformation.
Algorithms to accomplish this have recently been 
developed in the case of intersecting
hypersurfaces in a {\em complex\/} torus~\cite{CVVerschelde,Huber_Sturmfels}.  

This note is organized as follows:
In \S 2 we discuss some examples of fully real enumerative problems
for which multiplicity-free deformations play  a central role.
This technique is illustrated  in \S 3, where we show that there are
nine real  Veronese surfaces in ${\Bbb P}^5$ such that the $11010048$ 
planes meeting all nine are real.
We conclude with a discussion of the work of Fulton and of Ronga,
Tognoli, and Vust~\cite{Ronga_Tognoli_Vust} on the problem of conics
tangent to five conics and show that the multiplicities they introduce
are unavoidable.

\section{Effective Rational Equivalence} 

A common feature of many fully real enumerative problems is
multiplicity-free deformations of intersection cycles.
Effective rational equivalence is a precise formulation of this for
Grassmannians and flag varieties.

\subsection{Real effective rational equivalence}
Varieties will be quasi-projective, reduced, complex and defined 
over the real numbers, ${\Bbb R}$.
Let $X$ be a Grassmannian or flag variety, $G$ a linear algebraic
group which acts transitively on $X$, and 
$B$ a Borel subgroup of $G$.
The letters $U$ and  $V$ denote smooth rational varieties.
Let the real points $Y({\Bbb R})$ of a variety $Y$ be equipped with the
classical topology.

A subvariety $\Xi\subset U\times X$ (or $\Xi\rightarrow U$) with
generically reduced equidimensional fibres over $U$ is a {\em family of
(multiplicity-free) cycles on $X$ over $U$}.
We assume all families are $G$-stable; 
if $Y$ is a fibre of $\Xi$ over $U$, then so are all 
translates of $Y$.
Associating a point $u$ of $U$ to the fundamental cycle of the fibre
$\Xi_u$ determines a morphism
$\phi: U\rightarrow \Chow X$. 
Here, $\Chow X$ is the Chow variety of $X$ parameterizing cycles of the
same dimension and degree as $\Xi_u$ (\cite{Samuel}, \S I.9).
A priori, $\phi$ is only a function.
However, if $C\subset U$ is a smooth curve, then 
$\Xi|_C$ is flat and the 
canonical map of the Hilbert
scheme to the Chow variety (\cite{Mumford_Fogarty}, \S 5.4)
shows $\phi|_C$ is a morphism.
By Hartogs' Theorem on separate analyticity, $\phi$ is in fact a
morphism.
In fact, if $U$ is normal, then 
$\phi$ is a morphism~(\cite[\S 1]{Kollar_rational} 
or~\cite[\S 3]{Friedlander_Mazur}). 
For a discussion of Chow varieties in the analytic category (which
suffices for our purposes), see~\cite{Barlet}.

Any cycle $Y$ on $X$ is rationally equivalent to an integral linear
combination of Schubert classes.
As Hirschowitz~\cite{Hirschowitz} observed, this rational equivalence
occurs within the closure of $B\cdot Y$ in 
$\Chow X$ since $B$-stable cycles of $X$ ($B$-fixed points in
$B\cdot Y$) are integral linear
combinations of  Schubert varieties.
If any  coefficients in this linear combination exceed 1, this
stable cycle has multiplicities.

A family $\Xi\rightarrow U$ of multiplicity-free cycles on $X$ has {\em
effective rational equivalence} with {\em witness} $Z$ if there is
a cycle $Z\in \overline{\phi(U)}$ which is a sum of distinct Schubert
varieties, and hence multiplicity-free.
An effective rational equivalence is {\em real} if 
$Z\in \overline{\phi(U({\Bbb R}))}$ and each component of $Z$ is a
Schubert variety defined by a real flag.

Suppose  $\Xi_1\rightarrow U_1,\ldots,\Xi_b\rightarrow U_b$ are
$G$-stable families of multiplicity-free cycles on $X$.
By Kleiman's Transversality Theorem~\cite{Kleiman}, there is a
nonempty open set $U\subset \prod_{i=1}^b U_i$ consisting of
$b$-tuples $(u_1,\ldots,u_b)$ such that the fibres 
$(\Xi_1)_{u_1},\ldots,(\Xi_b)_{u_b}$ meet
generically transversally.
Let $\Xi\subset U\times X$ be the
resulting family of intersection cycles
and call $\Xi\rightarrow U$ the {\em intersection problem} given by 
$\Xi_1,\ldots,\Xi_b$.

\begin{thm}\label{thm:real_lines}\ 
Any intersection problem given by families of Schubert varieties 
in the Grassmannian of lines in projective space has real effective
rational equivalence.
\end{thm}

We present a synopsis of  the proof in~\cite{sottile_real_lines}:
Let $X$ be the Grassmannian of lines in ${\Bbb P}^n$ and suppose 
$\Xi\rightarrow U$ is an intersection problem given by families of
Schubert varieties.
A sequence $\Psi_0\rightarrow V_0,\ldots,\Psi_c\rightarrow V_c$ of
families of multiplicity-free cycles on $X$ is constructed with each
$V_i$ rational, where $\Psi_0\rightarrow V_0$ is the family
$\Xi\rightarrow U$, $V_c$ is a point, and $\Psi_c$ a union of
distinct real Schubert varieties.
For each $i=0,\ldots, c$, let ${\cal G}_i\subset \Chow X$ be
$\phi(V_i({\Bbb R}))$, the set of fibres of the
family $\Psi_i\rightarrow V_i$ over $V_i({\Bbb R})$.

Then ${\cal G}_i\subset \overline{{\cal G}_{i-1}}$:
For any $v\in V_i({\Bbb R})$ a family $\Gamma\rightarrow C$
of cycles is constructed with $C$ a smooth rational curve, the cycle
$(\Psi_i)_v$ a fibre over $C({\Bbb R})$, and all other fibres of
$\Gamma$ are fibres of $\Psi_{i-1}$.
This family induces a morphism $\phi:C\rightarrow \Chow X$,
which shows  
$(\Psi_i)_v\in  \overline{{\cal G}_{i-1}}$ since
$\phi(C({\Bbb R}))-\{(\Psi_i)_v\} \subset {\cal G}_{i-1}$.
It follows that 
$\Psi_c \in \overline{{\cal G}_0} = \overline{\phi(U({\Bbb R}))}$, 
showing $\Xi\rightarrow U$
has real effective rational equivalence.
\QED

It is typically difficult to describe an intersection of several
Schubert varieties.
While this task is easier when they are in special position, 
even this may be too hard.
It is better yet to consider the limiting position of intersection
cycles as the subvarieties being intersected degenerate to the point
of attaining excess intersection.
This is the aim of effective rational equivalence.
For example,  the `limit cycle' $\Psi_c$ of the previous proof
is generally not an intersection of Schubert varieties,
however, it is a deformation of such cycles.

An {\em enumerative problem} of degree $d$ is an intersection problem
$\Xi\rightarrow U$ with finite fibres of cardinality $d$.
It is {\em fully real\,} if there is a fibre
$\Xi_u$ with $u\in U({\Bbb R})$ consisting entirely of real points.
Here,  $u=(u_1,\ldots,u_b)$
with $u_i\in U_i({\Bbb R})$ and $\Xi_u$ is the transverse intersection
of the cycles $(\Xi_1)_{u_1},\ldots,(\Xi_b)_{u_b}$.

The set ${\cal M}\subset \mbox{\it Sym}^d X$ of degree $d$ zero cycles
consisting of $d$ distinct real points of $X$ is an open subset of
$(\mbox{\it Sym}^d X)({\Bbb R})$.
Thus $\Xi\rightarrow U$ is fully real if and only if it has real
effective rational equivalence.
Hence Theorem~\ref{thm:real_lines} has the following consequence:

\begin{cor}\label{cor:real_lines}
Any enumerative problem given by Schubert conditions on lines in
projective space is fully real.
\end{cor}

\subsection{Products in $A^*X$}
$X$ is the quotient $G/P$ of $G$ by a parabolic subgroup
$P$. 
A Schubert subvariety $\Omega_w\Fdot$ of $X$ is given by a complete
flag $\Fdot$ and a coset $w$ of the corresponding parabolic subgroup
of the symmetric group (\cite{Bourbaki_Groupes_IV}, Ch.~IV, \S 2.5).
Call $w$ the {\em type} of $\Omega_w\Fdot$.
A Schubert class $\sigma_w$ is the cycle class of $\Omega_w\Fdot$.

Let $\Xi_1\rightarrow U_1,\ldots,\Xi_b\rightarrow U_b$ be families of
cycles on $X$ giving an intersection problem $\Xi\rightarrow U$.
Then fibres of $\Xi\rightarrow U$ have cycle class $\prod_i \beta_i$,
where $\beta_i$ is the cycle class of fibres of 
$\Xi_i\rightarrow U_i$.
Suppose $\Xi\rightarrow U$ has effective rational equivalence with
witness $Z$. 
Let $c_w$ count the components of $Z$ of type $w$.
Since $Z$ is rationally equivalent to fibres of $\Xi\rightarrow U$, 
we deduce the formula in $A^*X$.
$$
\prod_{i=1}^b \beta_i\quad =\quad \sum_w c_w\cdot \sigma_w.
$$

\subsection{Pieri-type formulas}
Given a such product formula with each $c_w\leq 1$,
the action of a real Borel subgroup $B$ of $G$ shows
that the family of intersection cycles $\Xi\rightarrow U$ has real
effective rational equivalence:
Let $Y$ be a fibre of $\Xi\rightarrow U$ over a real point of $U$.
Then the closure of the orbit $B({\Bbb R})\cdot Y$ in 
$\Chow X({\Bbb R})$ contains a $B({\Bbb R})$-fixed point $Z$, as Borel's
fixed point Theorem (\cite{Borel_groups}, III.10.4), holds for 
$B({\Bbb R})$-stable real analytic sets.
Moreover, $Z$ is multiplicity-free as $c_w\leq 1$.

In  the Grassmannian of $k$-planes in 
${\Bbb P}^n$, a {\em special Schubert variety} is the locus of $k$-planes
having excess intersection with a fixed linear subspace.
A {\em special Schubert variety} of a flag variety is the inverse
image of a special Schubert variety in a Grassmannian projection.
Pieri's formula for
Grassmannians~(\cite{Griffiths_Harris},~\cite{Hodge_Pedoe}) and the
Pieri-type formulas for flag
varieties~(\cite{Lascoux_Schutzenberger_polynomes_schubert}%
,~\cite{sottile_pieri_schubert}) show that the coefficients $c_w$ in a
product of a Schubert class with a special Schubert class are either 0
or 1.
Thus any intersection problem given by a Schubert variety and a 
special Schubert variety has real effective rational equivalence.
We use this to prove the following theorem.

\begin{thm}\label{thm:three_special}
Any enumerative problem in any flag variety given by five Schubert
varieties, three of which are special, is fully real.
\end{thm}

{\sc Proof}.  
First pair each non-special Schubert variety with a
special Schubert variety.  
The associated families $\Xi\rightarrow U$ and $\Xi'\rightarrow U'$ of
intersection cycles have real effective rational equivalence with
witnesses $Z$ and $Z'$, respectively.

Since the coefficients $c_w$ in the Pieri-type formulas
are either 0 or 1, a zero-dimensional intersection of three real
Schubert varieties in general position where one is special is a single
real point. 
Considering components of $Z$ and $Z'$ separately, we see that if 
$Z$, $Z'$, and the third special Schubert variety $Y$ are in
general position with $Y$ real, then they intersect transversally with
all points of intersection real.
Suitably small deformations of $Z$ and $Z'$ into real fibres of $\Xi$
and $\Xi'$ preserve the number of real points of intersection, 
completing the proof.
\QED

\section{The Grassmannian of planes in ${\Bbb P}^5$}
The Grassmannian of planes in ${\Bbb P}^5$, ${\Bbb G}\,_{2,5}$, is
a 9-dimensional variety.
If $K$ is a plane in  ${\Bbb P}^5$, then the set  $\Omega(K)$ of
planes which meet $K$ is a hyperplane section of 
${\Bbb G}\,_{2,5}$ in its Pl{\"u}cker embedding.
Thus the number of planes which meet 9 general planes is the
degree of ${\Bbb G}\,_{2,5}$, which is 
$\frac{1!2!9!}{3!4!5!}=42$~\cite{Schubert_degree}.
This variety is the smallest dimensional
flag variety for which an analog of Corollary~\ref{cor:real_lines}
is not known.
We illustrate the methods of  \S 2 to prove the following result: 

\begin{thm}\label{thm:42_planes}
There are 9 real planes in ${\Bbb P}^5$ such that the 42
planes meeting all 9 are real.
\end{thm}

The Veronese surface in  ${\Bbb P}^5$  is the image of  ${\Bbb P}^2$
under the embedding induced by the complete linear system
$|{\cal O}(2)|$, and so it has degree 4.

\begin{cor}
There are 9 real Veronese surfaces in ${\Bbb P}^5$ such that the\
11010048 ($= 4^9\cdot 42$) planes meeting all 9 are real.
\end{cor}

{\sc Proof}. 
Let $x_{ij}$ , $1\leq i\leq j\leq 3$, be real coordinates
for ${\Bbb P}^5$.
For $t\neq 0$ 
\begin{multline}\label{eq:ideal_family}
\langle
\underline{x_{11}x_{33}}-t^4\,x_{13}^2,\  
\underline{x_{11}x_{22}}-t^2\,x_{12}^2,\  
\underline{x_{11}x_{23}}-t\, x_{12}x_{13},\\
\underline{x_{12}x_{33}}-t\, x_{13}x_{23},\  
\underline{x_{13}x_{22}}-t\, x_{12}x_{23},\  
\underline{x_{22}x_{33}}-t^2\,x_{23}^2 \rangle  
\end{multline}
generates the ideal of a Veronese surface, ${\cal V}(t)$
(cf.~\cite{Sturmfels_grobner_ULS}, p.~142), which is real for
$t\in{\Bbb R}$.
This family of Veronese surfaces is
induced by the (real) ${\Bbb C}^\times$-action
on the space of linear forms on ${\Bbb P}^5$:
$$
x_{ij}\quad \mapsto\quad t^{j-i}x_{ij}\qquad\mbox{for }t\in {\Bbb
C}^\times .
$$
The ideal of the special fibre ${\cal V}(0)$ of this family is generated
by the underlined terms, so ${\cal V}(0)$ is the union of the four
planes given by the ideals:
\begin{equation}\label{eq:four_planes}
\langle x_{11},x_{22},x_{33}\rangle \qquad
\langle x_{ii},x_{jj},x_{ij}\rangle,\quad  ij=12, 13, 23.
\end{equation}

By Theorem~\ref{thm:42_planes}, there exist 9 real planes
$K_1,\ldots,K_9$ such that $\bigcap_{i=1}^9 \Omega(K_i)$
is a transverse intersection consisting of 42 real planes.
This property of $K_1,\ldots, K_9$ is
preserved by small real deformations.
So for each $1\leq i\leq 9$, there is a 
neighborhood $W_i$ of $K_i$ in ${\Bbb G}\,_{2,5}({\Bbb R})$ such
that if $K'_i\in V_i$ for $1\leq i\leq 9$, then 
$\bigcap_{i=1}^9 \Omega(K'_i)$
is transverse and consists of 42 real planes.

For each $1\leq i\leq 9$, choose a set of  real coordinates for 
${\Bbb P}^5$ so that the four planes, $K_{i,j}$, for $j=1,2,3,4$,
defined by the ideals of~(\ref{eq:four_planes}) are in $W_i$.
In these same coordinates, consider the family ${\cal V}_i(t)$ of real
Veronese surfaces given by the ideals~(\ref{eq:ideal_family})
with special member 
${\cal V}_i(0) = K_{i,1} + K_{1,2}+ K_{i,3}+ K_{i,4}$.
If the sets of coordinates are chosen sufficiently generally, there
exists $\epsilon >0$ such that whenever $t\in (0,\epsilon)$, there are
exactly $4^9\cdot 42$ real planes meeting each of ${\cal V}_1(t),
\ldots,{\cal V}_9(t)$.

This is because there are  $4^9\cdot 42$ real planes meeting each
of ${\cal V}_1(0),\ldots,{\cal V}_9(0)$, as
$$
\bigcap_{i=1}^9 \left(\Omega(K_{i,1})+\Omega(K_{i,2})+
\Omega(K_{i,3})+\Omega(K_{i,4}) \right)
$$
is a transverse intersection consisting of $4^9\cdot 42$ real planes:
Since $K_{i,j}\in W_i$ for $1\leq i\leq 9$ and $1\leq j\leq 4$,
this follows if the $4^9$ sets of 42 planes 
$\bigcap_{i=1}^9\Omega(K_{i,l_i})$
given by all sequences $l_i$, where $1\leq l_i\leq 4$ for $1\leq i\leq 9$,
are pairwise disjoint. 
But this may be arranged when choosing the seta of coordinates.
\QED

\begin{lemma}\label{lemma:real_eff_rat_equiv}
The intersection problem of planes meeting 4 given planes in 
${\Bbb P}^5$ has real effective rational equivalence.
\end{lemma}

{\sc Proof of Theorem~\ref{thm:42_planes} using
Lemma~\ref{lemma:real_eff_rat_equiv}}.
Partition the 9 planes into two sets of 4, and a singleton.
Apply Lemma~\ref{lemma:real_eff_rat_equiv} to the intersection
problems $\Xi\rightarrow U$,  $\Xi'\rightarrow U'$ given by each 
set of 4, obtaining witnesses $Z$ and $Z'$.
Arguing as for Theorem~\ref{thm:three_special} completes
the proof.
\QED

{\sc Proof of Lemma~\ref{lemma:real_eff_rat_equiv}}.
We use an economical notation for Schubert varieties.
A partial flag $A_0\subsetneqq A_1\subsetneqq A_2\subset {\Bbb P}^5$
determines a Schubert subvariety of ${\Bbb G}\,_{2,5}$:
$$
\Omega(A_0,A_1,A_2) := \{H\in {\Bbb G}\,_{2,5}\,|\,
\dim H\bigcap A_i \geq i, \mbox{ for }i=0,1,2\}.
$$
If $A_j$ is a hyperplane in $A_{j+1}$ or if $A_j={\Bbb P}^5$, then 
it is no additional restriction for $\dim H\bigcap A_j\geq j$.
We omit such inessential conditions.
Thus, if $\mu\subsetneqq M\subsetneqq\Lambda$ is a partial flag, then 
$\Omega(\mu), \Omega(\dotc,M)$, and $\Omega(\mu,\dotc,\Lambda)$ are,
respectively, those planes $H$ which meet $\mu$,
those $H$ with $\dim H\cap M\geq 1$, and those $H\subset\Lambda$ which 
also meet $\mu$.

Let $\Xi\subset U\times {\Bbb G}\,_{2,5}$ be the intersection problem of
planes meeting four given planes.
Then $U\subset ({\Bbb G}\,_{2,5})^4$ is the set of 4-tuples of planes
$(K_1,K_2,K_3,K_4)$ such that $\bigcap_{i=1}^4 \Omega(K_i)$ is a
generically transverse intersection
and the fibre of $\Xi$ over $(K_1,K_2,K_3,K_4)$ is 
$\bigcap_{i=1}^4 \Omega(K_i)$.
We show $\Xi\rightarrow U$ has real effective rational equivalence by
exhibiting a family $\Psi\subset V\times {\Bbb G}\,_{2,5}$, 
satisfying the four conditions:
\begin{enumerate}
\item[(a)] $V$ is rational.  
In fact $V$ is a dense subset of Magyar's configuration variety
${\cal F}_D$~\cite{Magyar_Borel-Weil}, where $D$ is the diagram
$$
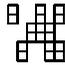
\DIAG
$$
\item[(b)] $V$ has a dense open subset $V^\circ$ such that the fibres of
$\Psi|_{V^{\circ}}$ are also fibres in the family $\Xi\rightarrow U$.
\item[(c)] $V$ has a rational subset $V'$ such that the fibres of
$\Psi|_{V'}$ are unions of distinct Schubert varieties, real for real
points of $V'$.
\item [(d)]$V'({\Bbb R})\subset\overline{V^\circ({\Bbb R})}$.
Hence $\phi(V'({\Bbb R}))\subset 
\overline{\phi(U({\Bbb R}))}$. 
Together with (c), this shows $\Xi\rightarrow U$ has 
effective rational equivalence.
\end{enumerate}

Let $V\subset ({\Bbb G}_{1,5})^3\times({\Bbb G}\,_{3,5})^3$ be the locus
of sextuples $(\mu_1,\mu_2,\lambda;M_1,M_2,L)$ such that 
$\mu_i\subset M_i$, $i=1,2$, $\mu_1,\mu_2\subset L$, 
$\lambda \subset M_1\bigcap M_2$, and 
$\mu_i\not\subset M_j, i\neq j$.
We illustrate the inclusions:
$$
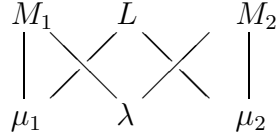
\begin{picture}(100,50)
\put(0,0){$\mu_1$}
\put(40,0){$\lambda$}
\put(85,0){$\mu_2$}
\put(0,38){$M_1$}
\put(40,38){$L$}
\put(85,38){$M_2$}
\put(5,10){\line(0,1){25}}
\put(90,10){\line(0,1){25}}
\put(40,10){\line(-1,1){25}}
\put(50,10){\line(1,1){25}}
\put(40,35){\line(-1,-1){11}}
\put(26,21){\line(-1,-1){11}}
\put(50,35){\line(1,-1){11}}
\put(64,21){\line(1,-1){11}}
\end{picture}
$$

Let $V^\circ\subset V$ be the dense locus where 
$\langle \mu_i,M_j\rangle ={\Bbb P}^5$
for $i\neq j$.
Then $\lambda = M_1\bigcap M_2$ and $L=\langle \mu_1,\mu_2\rangle$.
Let $V'\subset V$ be the locus where $\mu_1\bigcap \mu_2$ is a point,
so that $\langle M_1,M_2\rangle$ is a hyperplane.
Then $V'$ is rational and $V'({\Bbb R})\subset
\overline{V^\circ({\Bbb R})}$, proving (d).

We define the family $\Psi$.
For $v=(\mu_1,\mu_2,\lambda;M_1,M_2,L)\in V$, let $\Psi_v$ be
the cycle
\begin{multline}\label{eq:Psi_v}
\Omega(\mu_1)\bigcap\Omega(\dotc,M_2) \ +\ 
\Omega(\mu_2)\bigcap\Omega(\dotc,M_1) \ +\\ 
\{H\in \Omega(\mu_1,L)\,|\, H\bigcap \mu_2\neq \emptyset\} \ +\
\{H\in \Omega(\lambda,M_1)\,|\, \dim H\bigcap M_2\geq 1\}.
\end{multline}
Let $\Psi\subset {\Bbb G}\,_{2,5}\times V$ be the subvariety with  
fibre  $\Psi_v$ over points $v\in V$.

Suppose $v\in V^\circ$.
Since $L=\langle \mu_1,\mu_2\rangle$ and
$\mu_1\bigcap\mu_2=\emptyset$,
$$
\Omega(\mu_1)\bigcap\Omega(\mu_2) \ =\ 
\{H\in \Omega(\mu_1,L)\,|\, H\bigcap \mu_2\neq \emptyset\},
$$
as any plane meeting both $\mu_1$ and $\mu_2$ must intersect their
span $L$ in at least a line.
Similarly, $\Omega(\dotc,M_1)\bigcap\Omega(\dotc,M_2)$ is the fourth
term of the cycle (\ref{eq:Psi_v}): 
If $l_i$ is a line in $H\bigcap M_i$ for
$i=1,2$, then $l_1\bigcap l_2 \subset \lambda = M_1\bigcap M_2$.
Thus 
$\Psi_v=\bigcap_{i=1}^2\left(\Omega(\mu_i)+\Omega(\dotc,M_i)\right)$.
Since the pairs of subspaces $(\mu_1,\mu_2)$, $(M_1,M_2)$, and 
$(\mu_i,M_j)$ for $i\neq j$ are in general position, this intersection
is generically transverse.

We claim $\Psi_v$ is a fibre of $\Xi\rightarrow U$:
Let $K_i,K_i'$ for $i=1,2$ be planes such that 
$\mu_i=K_i\bigcap K'_i$ and $M_i=\langle K_i,K_i'\rangle$.
Then $\Omega(K_i)\bigcap\Omega(K'_i) = 
\Omega(\mu_i)+\Omega(\dotc,M_i)$:
If a plane $H$ meets
both $K_i$ and $K_i'$, either it meets their intersection $\mu_i$, or else
it intersects their span $M_i$ in at least a line.
Moreover, while $K_i,K_i'$ are not in general position, this
intersection is generically transverse as a proper intersection of a
Schubert variety with a special Schubert variety is necessarily
generically transverse (\cite{sottile_explicit_pieri}, \S 2.7).
Thus $\Psi_v = \Xi_{(K_1,K'_1,K_2,K'_2)}$, proving (b).

To show (c), let 
$v=(\mu_1,\mu_2,\lambda;M_1,M_2,L)\in V'$.
Set $p=\mu_1\bigcap \mu_2$, a point and 
$\Lambda=\langle M_1,M_2\rangle$, a hyperplane.
Then $\langle \mu_1,\mu_2\rangle$ is a plane $\nu$ contained in $L$
and $M_1\bigcap M_2$ is a plane $N$ containing $\lambda$.
We illustrate these inclusions:
$$
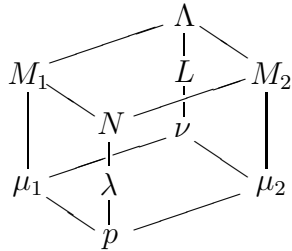
\begin{picture}(108,100)(-3,-5)
\put(33,0){$p$}
\put(-1,20){$\mu_1$}\put(32,18){$\lambda$}\put(91,20){$\mu_2$}
\put(31,41){$N$}\put(60,40){$\nu$}\put(-3,60){$M_1$}
\put(60,61){$L$}\put(89,60){$M_2$}\put(60,82){$\Lambda$}
\put(30,5){\line(-3,2){18}}\put(30,48){\line(-3,2){18}}
\put(88,25){\line(-3,2){20}}\put(88,70){\line(-3,2){18}}
\put(45,5){\line(3,1){40}}\put(45,48){\line(3,1){42}}
\put(13,25){\line(3,1){21}}\put(37,33){\line(3,1){20}}
\put(13,68){\line(3,1){43}}
\put(35.5,7){\line(0,1){9}}\put(35.5,29){\line(0,1){9}}
\put(64,48){\line(0,1){5}}\put(64,56){\line(0,1){3}}
\put(64,71){\line(0,1){8}}
\put(5,28){\line(0,1){29}}\put(95,28){\line(0,1){29}}
\end{picture}
$$
To complete the proof, we show $\Psi_v$ is the sum of 
Schubert varieties
\begin{multline*}
\Omega(\mu_1,\dotc,\Lambda) \ +\ \Omega(p,M_2) \ +\ 
\Omega(\mu_2,\dotc,\Lambda) \ +\ \Omega(p,M_1) \ +\\
\Omega(\dotc,\nu) \ +\ \Omega(p,L) \ +\ 
\Omega(\lambda,\dotc,\Lambda) \ +\ \Omega(\dotc,N).
\end{multline*}
First note that 
$$\Omega(\mu_1)\bigcap \Omega(\dotc,M_2)\ =\ 
\Omega(\mu_1,\dotc,\Lambda)  + \Omega(p,M_2):
$$
If $H\in \Omega(\mu_1)\bigcap \Omega(\dotc,M_2)$, then 
either $H\bigcap \mu_1\not\subset M_2$, so that 
$H\subset \langle \mu_1,M_2\rangle = \Lambda$, or else
$p\in H$ so that $H\in \Omega(p,M_2)$.
Similarly, we have 
$\Omega(\mu_2)\bigcap \Omega(\dotc,M_1)=
\Omega(\mu_2,\dotc,\Lambda)  + \Omega(p,M_1)$.
These intersections are generically transverse, as they are proper.

Furthermore, 
$$
\{H\in \Omega(\mu_1,L)\,|\, H\bigcap \mu_2\neq \emptyset\} 
\ =\  \Omega(\dotc,\nu) \ +\ \Omega(p,L):
$$
Either $H\bigcap \mu_1\bigcap \mu_2= \emptyset$, thus 
$\dim H\bigcap\langle \mu_1,\mu_2\rangle\geq 1$, and so
$H\in\Omega(\dotc,\nu)$, or else $p\in H$, so that 
$H\in\Omega(p,L)$.
Finally, 
$$
\{H\in \Omega(\lambda,M_1)\,|\, \dim H\bigcap M_2\geq 1\}\ =\ 
\Omega(\lambda,\dotc,\Lambda) \ +\ \Omega(\dotc,N):
$$
Either
$H\bigcap M_1\not\subset M_2$ thus  
$H\subset \langle M_1,M_2\rangle=\Lambda$ and so 
$H\in \Omega(\lambda,\dotc,\Lambda)$, or else
$\dim H\bigcap M_1\bigcap M_2\geq 1$, so that 
$H\in \Omega(\dotc,N)$.
\QED

Note that for $v\in V'$, the fibre $\Psi_v$ is {\em not} an
intersection of four Schubert varieties of type
$\Omega(K)$, for $K$ a plane:
The Schubert subvariety $\Omega(\dotc,N)$ is the locus of planes which
contain a line $l\subset N$ and hence it consists of all planes of the
form  $\langle q,l\rangle$, where $l\subset N$ is a line and 
$q\in {\Bbb P}^5\setminus l$ is a point.
Suppose $\Psi_v\subset \Omega(K)$ so that
$\Omega(\dotc,N)\subset\Omega(K)$. 
Then for every line $l\subset N$ and point $q\in {\Bbb P}^5\setminus l$,
we have $K\bigcap\langle q,l\rangle\neq\emptyset$.
This implies that  $K\bigcap l\neq \emptyset$ for every line 
$l\subset N$, and hence that 
 $\dim K\bigcap N\geq 1$.
Similarly, $\dim K\bigcap \nu\geq 1$, and so 
$K\bigcap N\bigcap\nu\neq\emptyset$, thus $p\in K$.
This shows $\Omega(p)\subset \Omega(K)$
and so if 
$\Psi_v\subset \Omega(K_1)\bigcap\Omega(K_2)
\bigcap\Omega(K_3)\bigcap\Omega(K_4)$,
then this intersection must contain $\Omega(p)$.
Hence  $\Psi_v$ is a proper subset of the intersection.

If $a_i=\dim A_i$, then $\sigma_{a_1a_2a_3}$ is the rational
equivalence class of $\Omega(A_0,A_1,A_2)$.
By the observation of  \S 2.2, Lemma~\ref{lemma:real_eff_rat_equiv}
implies the formula in $A^*{\Bbb G}\,_{2,5}$:
$$
(\sigma_{245})^4\ =\ 3\cdot \sigma_{035}\ +\ 
2\cdot \sigma_{125}\ +\ 3\cdot \sigma_{134},
$$
which may be determined by other means from the classical Schubert calculus.

\section{Real Plane Conics}

In 1864 Chasles~\cite{Chasles} showed there are 3264 plane conics
tangent to five general conics.
Fulton~\cite{Fulton_introduction_intersection} asked how many of the
3264 conics tangent to five general (real) conics can be real.
He later determined that all can be real, but did not publish that
result. 
More recently, Ronga, Tognoli, and Vust~\cite{Ronga_Tognoli_Vust}
rediscovered this result.
We conclude this note with an outline of their work.
The author is grateful to Bill Fulton and Felice Ronga for
explaining these ideas.

Let $X$ be the variety of complete plane conics, a smooth variety of
dimension 5.
Let the hypersurfaces $H_p$, $H_l$, and $H_C$, be, respectively those
conics containing a point $p$, those tangent to a line $l$, and those
tangent to a conic, $C$.
If $\hat{p},\hat{l}$, and $\widehat{C}$ are, respectively, their cycle
classes in $A^1 X$, then 
$$
\widehat{C}\ =\ 2 \hat{p}\ +\ 2 \hat{l},
$$
which may be seen by degenerating a conic into two lines.
Then the number of conics tangent to five general conics is the degree
of 
$$
\widehat{C}^5 \ =\ 32 (\hat{p}\,^5\ +\ 5\hat{p}\,^4\cdot\hat{l}\ +\ 
10\hat{p}\,^3\cdot\hat{l}\,^2\ +\ 10\hat{p}\,^2\cdot\hat{l}\,^3\ 
+\ 5\hat{p}\cdot\hat{l}^4\ +\ \hat{l}\,^5).
$$
The monomials $\hat{p}\,^j\cdot\hat{l}^{5-j}$ for $j=0,\ldots,5$, 
have degrees $1,2,4,4,2,1$, giving Chasles' number of 
$32(1+10+40+40+10+1)=3264$~\cite[\S 9]{Kleiman_pspum}.

\begin{thm}[Fulton, Ronga-Tognoli-Vust]
There are five real conics in general position such that all of the 3264
conics tangent to the five conics are real.
\end{thm}

{\sc Proof.} 
The strategy is to realize the five conics as  a deformation of
five degenerate conics giving a maximal number of
real conics.
The first step is to show that for each $j$, there are $j$ lines and
$5-j$ points such that the  $2^{\min\{j, 5-j\}}$ conics tangent to the
lines and containing the points are real.
In~\cite{Ronga_Tognoli_Vust}, this step is done explicitly with a precise
determination of which configurations of points and lines are
`maximal'; that is, have all solutions real.
Remarkably, there are five lines $l_1,\ldots,l_5$ and five real points
$p_1,\ldots,p_5$ with $p_i\in l_i$ such that each of the 32 terms in
$$  
\bigcap_{i=1}^5 \left (H_{p_i}\ +\ H_{l_i} \right)
$$	
is a transverse intersection with all points of intersection real.
Such a configuration is illustrated in Figure~\ref{fig:one}.
\setcounter{figure}{0}
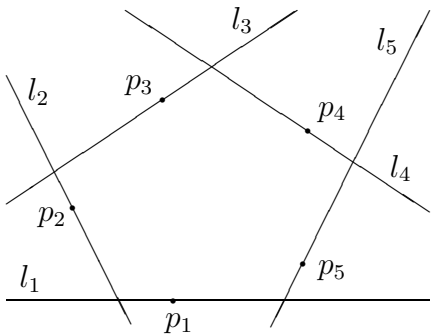
\begin{figure}[htb]
$$
\begin{picture}(170,130)(-5,-5)
\put(60,0){$p_1$}\put(63,10){\circle*{2}}
\put(5,15){$l_1$}
\put(0,10){\line(1,0){160}}

\put(12,40){$p_2$}\put(8,86){$l_2$}\put(25,45){\circle*{2}}
\put(25,45){\line(-1,2){25}}\put(25,45){\line(1,-2){22}}

\put(45,90){$p_3$}\put(85,113){$l_3$}\put(59,86){\circle*{2}}
\put(59,86){\line(-3,-2){59}}\put(59,86){\line(3,2){51}}

\put(118,80){$p_4$}\put(145,57){$l_4$}\put(114,74){\circle*{2}}
\put(114,74){\line(3,-2){46}}\put(114,74){\line(-3,2){69}}

\put(118,20){$p_5$}\put(140,105){$l_5$}\put(112,24){\circle*{2}}
\put(112,24){\line(1,2){48}}\put(112,24){\line(-1,-2){12}}
\end{picture}
$$
		\caption{A Maximal Configuration}\label{fig:one}
\end{figure}

The maximality of such a configuration is stable under small real
deformations of its points and lines. 
Thus we may choose real lines $l_1',\ldots,l_5'$ where
\begin{enumerate}
\item $p_i\in l_i'$ and $l_i'$ is distinct from $l_i$, for $i=1,\ldots,5$,
\item Any configuration obtained from a maximal configuration by
substituting some primed lines for the corresponding unprimed lines is
maximal.
\item The lines $l_i$ and $l_i'$ partition the real tangent directions
at $p_i$ into two intervals.
The configurations described in condition (2) give finitely many
real conics passing through $p_i$.
We require that all tangent directions to these conics at $p_i$ 
lie within the interior of {\em one} of these two intervals.
\end{enumerate}

The relation $\widehat{C} = 2\hat{p} + 2\hat{l}$ may be obtained by
considering a conic $C$ near a degenerate conic consisting of
two lines $l,l'$ meeting at a point $p$, and a pencil of
conics.\footnote{This version of this manuscript does not contain all of
the postscript files of the original, in particular, it does not have an
illustration of this degenerate conics $C$ and the nearby conics.
To visualize this, think of a real hyperbola in ${\Bbb R}^2$
defined by $x^2 -y^2=t$, where $t$ is a small positive real number.
Then the two lines are $x=\pm y$ and the point $p$ is the origin.
The condition that the real tangent line to $Q$ at $p$ does not
intersect $C$ means the absolute value of its slope exceeds 1.}
For any conic $q$ in that pencil tangent to one of the lines, there is
a nearby conic $q'$ in that pencil tangent to $C$.
However, for every conic $Q$ in the pencil containing $p$, there are
{\em two} nearby conics $Q', Q''$  in that pencil tangent to $C$.
Moreover, if $Q$ is real, then $Q'$ and $Q''$ are real if and only if
the real tangent line to $Q$ at $p$ does not intersect $C$.

By condition (3), we may choose real conics $C_1,\ldots,C_5$ with $C_i$
near the degenerate conic $l_i + l_i'$ and, if
$Q$ is a conic in 
\begin{equation}\label{eq:last}
\bigcap_{i=1}^5 \left( H_{p_i}\ +\ H_{l_i}\ +\ H_{l_i'}\right)
\end{equation}
containing $p_i$, the the real tangent line to $Q$ at $p_i$ does not
intersect $C_i$.
If, in addition, the conics $C_i$ are sufficiently close to each
degenerate conic, then there will be 3264 real conics tangent to each
of $C_1,\ldots,C_5$.

Indeed, suppose $H_{C_1}$ replaces $H_{p_1} + H_{l_1} + H_{l_1'}$ in
the intersection (\ref{eq:last}).
Then for any conic $q$ in (\ref{eq:last}) that is tangent to either
$l_1$ or $l_1'$, there is a nearby real conic $q'$ tangent to $C$ 
which satisfies the other conditions on $q$
(since these other conditions determine a pencil of conics).
If $Q$ is a conic in (\ref{eq:last}) containing $p_1$, 
then there are two nearby real conics $Q'$ and
$Q''$  tangent to $C$ which satisfy the other conditions on $Q$.
Similarly, if $H_{C_2}$ now replaces 
$H_{p_2} + H_{l_2}+ H_{l_2'}$ in the new intersection
$H_{C_1}\cap\bigcap_{i=2}^5\left( H_{p_i}+H_{l_i}+H_{l_i'}\right)$, then
each conic tangent to $l_2$ and $l_2'$ gives a conic tangent to
$C_2$, but each conic through $p_2$ gives two conics tangent to $C_2$.
Replacing $H_{C_3}, H_{C_4}$, and $H_{C_5}$ in turn completes the
argument.
\QED

This proof used the effective rational equivalence:
$$
H_C \ \sim\ 2 H_p\ +\ H_l\ +\ H_{l'},
$$
where $l, l'$ form a degenerate conic with $p= l\bigcap l'$.
This deformation to a cycle having multiplicities (the
coefficient 2 of $H_p$) is unavoidable:
The variety $X$ and thus $\Chow X$ has an action of
$G=PGL(3,{\Bbb C}\,)$.
The locus of hypersurfaces $H_C$ on $\Chow X$ is a single
5-dimensional $G$-orbit.
This family cannot have effective rational equivalence.
If $Z$ is a cycle in the closure of this locus, then $Z$ is in a
$G$-orbit of dimension at most 4.
Thus if $Z= H_p+H_{p'}+H_l+H_{l'}$, then the dimension of the
$G$-orbit of $(p,p',l,l')$ in the product of ${\Bbb P}^2$'s and their
duals is at most 4.
But this is impossible unless either $p=p'$ or $l=l'$.

\end{document}